\documentclass[reqno]{amsart}

\usepackage{amsmath}
\usepackage{amsthm}
\usepackage{amssymb}
\usepackage{mathrsfs} 
\usepackage{mathtools}

\usepackage{bm}
\usepackage{courier}
\usepackage{fnpct}
\usepackage{stmaryrd}
\usepackage[dvipsnames]{xcolor}
\usepackage{xcolor}

\definecolor{CBred}{RGB}{228,26,28}
\definecolor{CBblue}{RGB}{55,126,184}
\definecolor{CBgreen}{RGB}{77,175,74}
\definecolor{CBpurple}{RGB}{152,78,163}
\definecolor{CBorange}{RGB}{255,127,0}
\definecolor{CByellow}{RGB}{255,255,51}

\usepackage{subcaption}
\usepackage{stmaryrd}
\usepackage{hyperref}
    \hypersetup{
        colorlinks,
        linkcolor={blue!50!white},
        citecolor={blue!50!white},
        urlcolor={blue!80!white}
    }
\usepackage{comment}

\usepackage{tikz}
\usepackage{pgfplots}
\usetikzlibrary{arrows.meta}
\usetikzlibrary{shapes.misc}
\usepackage{ifthen}
\usepackage{tkz-graph}
\usepackage{subcaption}
\usepackage{caption}
\usetikzlibrary{positioning}
\usetikzlibrary{patterns}
\usetikzlibrary{calc}

\pgfdeclarelayer{background}
\pgfdeclarelayer{foreground}
\pgfsetlayers{background,main,foreground}

\numberwithin{equation}{section}

\newcommand{\Z}{\mathbb{Z}}

\renewcommand{\H}{\mathcal{H}}

\newcommand{\diam}{\mathrm{diam}}

\newcommand{\supp}{\mathrm{supp}}

\newcommand{\floor}[1]{\left\lfloor #1\right\rfloor}
\newcommand{\ceil}[1]{\left\lceil #1\right\rceil}

\newcommand{\restrict}{\upharpoonright}

\newcommand{\spn}{\textrm{span}}

\newcommand{\st}{:}

\newtheorem{theorem}{Theorem}[section]
\newtheorem{lemma}[theorem]{Lemma}
\newtheorem{corollary}[theorem]{Corollary}
\newtheorem{proposition}[theorem]{Proposition}

\theoremstyle{definition}
\newtheorem{question}[theorem]{Question}
\newtheorem{openproblem}[theorem]{Open Problem}

\newtheorem{remark}[theorem]{Remark}

\newtheorem{definition}[theorem]{Definition}
\newtheorem{example}[theorem]{Example}

\author[B. ~Cody]{Brent Cody}
\address[Brent Cody]{ 
Virginia Commonwealth University,
Department of Mathematics and Applied Mathematics,
1015 Floyd Avenue, PO Box 842014, Richmond, Virginia 23284, United States
} 
\email[B. ~Cody]{bmcody@vcu.edu} 
\urladdr{http://www.people.vcu.edu/~bmcody/}

\author[R. ~Detore]{Rose Detore}
\address[Rose Detore]{ 
Virginia Commonwealth University,
Department of Mathematics and Applied Mathematics,
1015 Floyd Avenue, PO Box 842014, Richmond, Virginia 23284, United States
} 
\email[R. ~Detore]{detorer@vcu.edu} 
\urladdr{}

\date{\today}

\pgfplotsset{compat=1.18}
\begin{document}

\title[Metric general position invariants and perfection]{Metric general position extensions of classical graph invariants and perfection}

\thanks{}

\begin{abstract}

We introduce a two-parameter framework that refines several classical graph invariants by imposing higher-order constraints along bounded-length geodesics. For integers \(k,d\ge1\), a vertex set is called \(k,d\)-independent if every shortest path of length at most \(d\) contains fewer than \(k\) vertices of the set, giving rise to corresponding \(k,d\)-generalizations of independence, chromatic, clique, and domination invariants. We establish basic bounds and monotonicity prop for these new invariants. In particular, we observe that a lower bound for the $k,d$-chromatic number can be obtained in terms of the size of a largest $k,d$-clique. We say that a graph $G$ is $k,d$-perfect if for every induced subgraph $H$, the $k,d$-clique lower bound of $H$ equals the $k,d$-chromatic number of $H$. We show that, although the $2,d$-invariants reduce to classical invariants in the $d$-th graph power, the theory of $2,d$-perfect graphs does not reduce to classical perfection in graph powers. Exact formulas are obtained for the \(k,d\)-chromatic number and $k,d$-clique number of paths and cycles, which leads to a complete characterization of which paths and cycles are $k,d$-perfect. 



\end{abstract}



\maketitle


\section{Introduction and definitions of $k,d$-invariants}\label{section_introduction}

Many fundamental notions in graph theory, such as independent sets, colorings, and cliques, are defined by imposing constraints on pairs of vertices. In this paper, we place these classical concepts within a broader two-parameter framework that incorporates both higher-order interactions and metric structure. The first parameter, $d$, is used to impose a geodesic condition by considering shortest paths of bounded length. This connects our work to a well-established range of distance-based notions, including distance-$d$ independent sets \cite{MR2378044}, packing colorings \cite{MR4126189}, general $d$-position sets \cite{MR4341189}, and related metric separation conditions, and is naturally linked to the theory of graph powers \cite{MR4688000}. The second parameter, $k$, shifts the focus from pairwise restrictions to constraints on $k$-element vertex sets. This direction extends the recent and rapidly growing literature on the general position problem for graphs \cite{MR4281067, MR4265041, MR4154901, MR4019752, MR3849577, MR3879620} by effectively replacing the collinearity of triples with constraints on larger vertex sets.

Taken together, these two parameters lead to the notion of $k$-general $d$-position sets \cite{MR4854543}, equivalent to the \emph{$k,d$-independent sets} (see Definition \ref{definition_main}) studied in this work. These sets simultaneously generalize classical independent sets and general position sets by controlling higher-order geodesic structure rather than simple pairwise distances. This unified perspective suggests the existence of a broader family of $k,d$-versions of classical graph invariants, including chromatic, clique, and domination numbers, that interact with $k,d$-independence in much the same way their classical counterparts interact with standard independence. To date, however, this perspective has appeared only in a limited form, most notably in the notion of general position colorings \cite{Chandran25}, defined for $k=3$ and $d=\mathrm{diam}(G)$.

In the present work, we introduce a comprehensive framework that places these ideas into a coherent theory. We make these interactions explicit and reveal new behavior even in basic graph classes such as paths and cycles. To begin our investigation, we define the $k,d$-analogues of four primary graph invariants.

\begin{definition}\label{definition_main}
Let $k$ and $d$ be positive integers and let $G$ be a graph.
\begin{enumerate}
    \item A set $S\subseteq V(G)$ is \emph{$k,d$-independent in $G$} if every shortest path in $G$ of length at most $d$ contains fewer than $k$ vertices of $S$.  
    The \emph{$k,d$-independence number of $G$} is
    \[
    \alpha^k_d(G)
    =\sup\left\{|S| \,\middle|\, \text{$S$ is a $k,d$-independent set in $G$}\right\}.
    \]

    \item For a positive integer $c$, a vertex coloring $f\colon V(G)\to[c]$ is a \emph{$k,d$-proper $c$-coloring of $G$}, or just a \emph{$k,d$-proper coloring}, if every color class of $f$ is $k,d$-independent; equivalently, every shortest path of length at most $d$ in $G$ contains fewer than $k$ vertices of the same color.  
    The \emph{$k,d$-chromatic number of $G$} is
    \[
    \chi^k_d(G)=\min\left\{c \,\middle|\, \text{$G$ admits a $k,d$-proper $c$-coloring}\right\}.
    \]

    \item A set $S\subseteq V(G)$ is a \emph{$k,d$-clique in $G$} if every $k$-element subset of $S$ is contained in some shortest path of $G$ of length at most $d$.  
    Define
    \[
    \omega^k_d(G)=\sup\left\{|S| \,\middle|\, \text{$S$ is a $k,d$-clique in $G$}\right\}.
    \]

    \item A set $D\subseteq V(G)$ is a \emph{$k,d$-dominating set in $G$} if for every vertex $v\in V(G)\setminus D$ there exists a shortest path in $G$ of length at most $d$ containing $v$ and satisfying $|V(P)\cap D|\ge k-1$.  
    The \emph{$k,d$-domination number of $G$} is
    \[
    \gamma^k_d(G)=\min\left\{|D| \,\middle|\, \text{$D$ is a $k,d$-dominating set in $G$}\right\}.
    \]
\end{enumerate}
\end{definition}

The parameters defined above extend the classical graph invariants. In the case $(k,d)=(2,1)$, we recover the standard parameters $\alpha^2_1(G)=\alpha(G)$, $\chi^2_1(G)=\chi(G)$, $\omega^2_1(G)=\omega(G)$, and $\gamma^2_1(G)=\gamma(G)$. For $k=2$ and $d \ge 1$, these invariants align with the metric structure of graph powers. Recall that for an integer $d\geq 1$, the \emph{$\ell$-th power of $G$}, denoted by $G^\ell$ is a graph with vertex set $V(G)$ such that $uv\in E(G^\ell)$ if and only if $u$ and $v$ lie on a shortest path of length at most $\ell$. For every $d \ge 1$, we have $\alpha^2_d(G)=\alpha(G^d)$, $\chi^2_d(G)=\chi(G^d)$, $\omega^2_d(G)=\omega(G^d)$, and $\gamma^2_d(G)=\gamma(G^d)$, where $G^d$ denotes the \emph{$d$-th power of $G$}. Just as in the classical setting, a simple counting argument shows that $\chi^k_d(G)\alpha^k_d(G) \ge |V(G)|$ for every graph $G$ (see Proposition \ref{proposition_basic}).

\begin{remark}\label{remark_general_framework}
The definition of $k,d$-independence can be viewed as a specific instance of a more general framework. Let $\mathcal{F}$ be an arbitrary family of subsets of $V(G)$. We may say that a set $S \subseteq V(G)$ is \emph{$k,\mathcal{F}$-independent} if for every $A \in \mathcal{F}$ we have $|S \cap A| < k$. Under this framework, our $k,d$-independence corresponds to the case where $\mathcal{F}$ is the set of all subsets of vertices that lie on a shortest path of length at most $d$. While some of our structural results could be extended to this generalized context (essentially treating the problem as independence in hypergraphs), the primary focus of this article is to exploit the specific metric structure of $G$. Thus, we retain the concrete definition above.
\end{remark}

For a set $S$ to be $k,d$-independent, vertices are permitted to be nearby provided they do not ``congest'' any shortest path of length at most $d$.  By extending this path-occupancy constraint to colorings, cliques, and domination, we recover a structural harmony that mirrors classical theory while also exhibiting novel behaviors. Although $k,d$-independent sets coincide with the \emph{$k$-general $d$-position sets} of \cite{MR4854543}, we adopt the terminology of $k,d$-independence to emphasize their functional role as higher-order analogues of independent sets within this broader framework.

Certain extreme parameter regimes will often be excluded from our results. For example, if $k > d+1$, then every shortest path of length at most $d$ contains fewer than $k$ vertices (by cardinality considerations), and consequently, every vertex set is $k,d$-independent, yielding the trivial values $\alpha^k_d(G)=|V(G)|$, $\chi^k_d(G)=1$, $\omega^k_d(G)=\min\{k-1,|V(G)|\}$, and $\gamma^k_d(G)=|V(G)|$. For this reason, many of our results begin by assuming $k\leq d+1$.

\begin{remark}\label{remark_hypergraph}
Let us note that the $k,d$-invariants defined in this article are equivalent to certan hypergraph invariants in an appropriate setting. Given a graph $G$, we can define a $k$-uniform hypergraph $\H_{k,d}(G)$ on $V(G)$ where $e\subseteq V(G)$ is a hyperedge if and only if $|e|=k$ and $e$ is contained in some shortest path of length at most $d$ in $G$. Then, for example, a set $S\subseteq V(G)$ is $k,d$-independent in $G$ if and only if it contains no hyperedge of $\H_{k,d}(G)$. The $k,d$-invariants defined above correspond to classical hypergraph invariants of $\H_{k,d}(G)$. Indeed, as shown in \cite{CodyDetore_powers}, the assignment $G \mapsto \mathcal{H}_{k,d}(G)$ defines a functor from the category of graphs with isometric embeddings to the category of $k$-uniform hypergraphs. However, let us offer two reasons for why the perspective and results presented in the current article are new, and do not simply reduce to facts about hypergraphs.
\begin{enumerate}
\item\label{item_justification_metric_structure}  In subsequent sections, we exploit the metric structure of $G$ to derive bounds and exact values that lie beyond the scope of general hypergraph theory. For example, the hypergraph greedy bound, $\chi^k_d(G)\le \Delta(\H_{k,d}(G))+1$, is a rather trivial and standard bound derived in terms of the maximum degree of $\H_{k,d}(G)$ (see Corollary \ref{corollary_greedy_hypergraph_bound}). Whereas the metric greedy bound, $\chi^k_d(G)\le \Delta(G^{\,d-k+2})+1$, is obtained by expoiting the metric structure of $G$ (see Proposition \ref{proposition_greedy_metric_bound}). We prove that these bounds are mutually incomparable by considering some explicit examples. Thus, the theory presented here does not reduce to general hypergraph theory.
\item While some of the general results in Section \ref{section_simple_bounds_and_kd_perfect_graphs} may be viewed as consequences of folklore results on the associated hypergraph $\H_{k,d}(G)$, and do not depend on additional properties of the underlying geodesic structure of $G$, related statements in the existing literature are dispersed and often formulated using differing conventions. To ensure clarity and self-containment, we therefore include explicit proofs of these results. 
\end{enumerate}
\end{remark}

In Section \ref{section_basic_properties}, we establish foundational properties of $k,d$-sets and invariants, including domination-type inequalities such as
\[
\omega^k_d(G)\le |V(G)|-\gamma^k_d(G)+(k-1),
\]
as well as results relating $k,d$-invariants of a graph to those of its isometric subgraphs, and results on the monotonicity of the invariants with respect to the parameters $k$ and $d$. Generalizing the clique lower bound on the chromatic number, as an easy consequence of the definitions, we obtain the $k,d$-clique lower bound (Proposition \ref{proposition_basic}):
\begin{equation}
\left\lceil\frac{\omega^k_d(G)}{k-1}\right\rceil \le \chi^k_d(G). \label{equation_clique_lower_bound_intro}
\end{equation}

In Section \ref{section_kd_clique_lower_bound}, motivated by the classical theory of perfect graphs \cite{MR2233847, MR309780}, we introduce \emph{$k,d$-perfect graphs} as those for which equality holds in (\ref{equation_clique_lower_bound_intro}) for every induced subgraph (see Definition \ref{definition_kd_perfect}). It easily follows from the definitions that a graph is $2,1$-perfect if and only if it is perfect, and thus this notion is a direct generalization. We prove that if $G^d$ is perfect then $G$ must be $2,d$-perfect (Proposition \ref{proposition_G^d_perfect_implies_G_is_2d_perfect}), which allows us to easily conclude via classical results that every strongly chordal graph is $2,d$-perfect for all $d \ge 1$ (Corollary \ref{corollary_strongly_chordal}), while chordal graphs are $2,d$-perfect for all odd $d \ge 1$ (Corollary \ref{corollar_chordal}). However, the converse is false: there is a graph $G$ which is $2,2$-perfect such that $G^2$ is not perfect (Proposition \ref{proposition_counter_example}). Hence, $2,d$-perfectness does not reduce to perfection of graph powers. We prove that every finite tree is $3,2$-perfect (Proposition \ref{proposition_trees_are_32_perfect}), a result that does not follow from classical perfection and whose proof underscores the higher-order and metric nature of our framework.


The remainder of the paper concers $k,d$-perfection of paths and cycles. In Section \ref{section_coloring paths}, we focus on paths. We derive an exact formula for the $k,d$-chromatic number of the finite path $P_n$, as well as for the two-way infinite path $P_\infty$. Similarly, we easily obtain formulas for $\omega^k_d(P_n)$ and $\omega^k_d(P_\infty)$. These formulas lead immediately to the conclusion that all paths are $k,d$-perfect for every admissible choice of $k$ and $d$. Thus, in the path setting, the higher-order and metric constraints interact in a remarkably rigid way, and the resulting behavior is realtively easily understood.

Section \ref{section_coloring_cycles} undertakes a parallel but substantially more delicate analysis for cycles. We obtain an explicit formula for $\chi^k_d(C_n)$ which generalizes a result of Kramer \cite[Theorem 1]{MR313102} on the distance-$d$ chromatic number of $C_n$. Our proof uses \emph{maximally even} sets of vertices in $C_n$, which originated in Clough and Douthett's work in music theory \cite{MR4875435, CloughDouthett}. Also in Section \ref{section_coloring_cycles}, we give an explicit formula for $\omega^k_d(C_n)$, and as an application of these formulas, we prove a complete characterization of which cycles are $k,d$-perfect. In the case $k=2$, the resulting picture closely mirrors the classical theory: for fixed $d$, small cycles are $2,d$-perfect, while larger cycles are $2,d$-perfect precisely when their lengths satisfy a simple divisibility  (Theorem \ref{theorem_2d_perfect_cycles}). For $k\ge3$, however, the situation changes (see Theorem \ref{theorem_kd_perfect_cycles}). For example, while all sufficiently short cycles are $k,d$-perfect, the behavior for larger cycles depends in a subtle way on the arithmetic relationship between $k$ and $d$. In particular, when $k-1$ divides $d+1$ there are infinitely many cycles that fail to be $k,d$-perfect, whereas when this divisibility fails, only finitely many cycles are not $k,d$-perfect. This leads to several phenomena that run counter to classical intuition. For example, as the reader may recall, for $n\geq 4$, $C_n$ is perfect if and only if $n$ is even. Whereas for $n\geq 3$, \emph{all} cycles $C_n$ are $3,2$-perfect. As another example, Corollary \ref{corollary_3d_perfect_cycles} reveals that for $n\geq 3$ the following hold.

\begin{itemize}
    \item $C_n$ is $3,2$-perfect.
    \item $C_n$ is $3,3$-perfect if and only if $n< 7$ or ($n\geq 7$ and $2\mid n$).
    \item $C_n$ is $3,4$-perfect if and only if $n\notin\{7\}$.
    \item $C_n$ is $3,5$-perfect if and only if ($n < 10$ and $n\notin\{7\}$) or ($n\geq 10$ and $3\mid n$).
    \item $C_n$ is $3,6$-perfect if and only if $n\notin\{7,10,11,13,17\}$.
    \item $C_n$ is $3,7$-perfect if and only if ($n<14$ and $n\notin\{7,10,11,13\}$) or ($n\geq 14$ and $4\mid n$).
\end{itemize}

Our results reveal an intricate structure which points to many natural structural and extremal questions that remain unresolved even for many basic families of graphs. We conclude the paper in Section~\ref{section_questions} with several open problems motivated by the phenomena uncovered throughout the article.

\section{Inequalities, isometric subgraphs and monotonicity}\label{section_simple_bounds_and_kd_perfect_graphs}

\subsection{Basic properties of $k,d$-sets and invariants}\label{section_basic_properties}

Let us begin with a few straightforward results that follow easily by extending the corresponding classical arguments.

\begin{proposition}\label{proposition_maximal_gp_implies_domination}
    Suppose $2\leq k\leq d+1$ are integers. Every maximal $k,d$-independent set in $G$ is a $k,d$-dominating set. Hence, 
    \[\gamma^k_d(G)\leq \alpha^k_d(G).\]
\end{proposition}

\begin{proof}
Let $S$ be a maximal $k,d$-independent set in $G$. For any $v\in V(G)\setminus S$ the set $S\cup\{v\}$ is not $k,d$-independent, and so there is a shortest path $P$ of length at most $d$ containing $v$ and containing $k-1$ vertices of $S$. Thus, $S$ is $k,d$-dominating.
\end{proof}

\begin{proposition}\label{proposition_clique_domination_bound}
    Suppose $2\leq k\leq d+1$ are integers. For every graph $G$,
    \[\omega^k_d(G)\leq |V(G)|-\gamma^k_d(G)+(k-1).\]
\end{proposition}

\begin{proof}
Let $Q$ be a $k,d$-clique with $|Q|=\omega^k_d(G)$. If $|Q|<k-1$ then the result follows trivially. So, suppose $|Q|\geq k-1$. Fix any set $A\subseteq Q$ with $|A|=k-1$ and let $D=A\cup(V(G)\setminus Q)$. Then $D$ is a $k,d$-dominating set because if $v\notin D$ then $v\in Q\setminus A$, so the set $A\cup\{v\}$, being a $k$-subset of $Q$, must be contained in a shortest path of $G$ with length at most $d$. Thus, $\gamma^k_d(G)\leq|D|=|V(G)|-|Q|+(k-1)$, which implies $\omega^k_d(G)\leq |V(G)|-\gamma^k_d(G)+(k-1)$.
\end{proof}

\begin{proposition}\label{proposition_basic}
For every graph $G$ and integers $k,d\geq 1$,
\[\ceil{\frac{|V(G)|}{\alpha^k_d(G)}}\leq \chi^k_d(G)\leq \ceil{\frac{|V(G)|-\alpha^k_d(G)}{k-1}}+1\]
and
\begin{align}
    \ceil{\frac{\omega^k_d(G)}{k-1}}\leq \chi^k_d(G).\label{equation_clique_bound_on_chromatic}
\end{align}
\end{proposition}

\begin{proof}
For the first inequality, let $f: V(G)\to [c]$ be a $k,d$-proper coloring $c$-coloring where $c=\chi^k_d(G)$. Then the color classes form a partition of $V(G)$, so
\[|V(G)|=\sum\limits_{i=1}^c|f^{-1}(i)|\leq\sum\limits_{i=1}^c\alpha^k_d(G)=\alpha^k_d(G)\chi^k_d(G).\]

For the second inequality, let $c=\ceil{\frac{|V(G)|-\alpha^k_d(G)}{k-1}}$. To see that $\chi^k_d(G)\leq c+1$, assign one color to the vertices of a particular $k,d$-independent set $S$ of size $\alpha^k_d(G)$ in $V(G)$. Then, partition $V(G)\setminus S$ into $c$ blocks, each of size $k-1$, plus potentially one additional block of size less than $k-1$. Each block is $k,d$-independent and is assigned a new color. In total, this coloring uses $c+1$ colors, hence $\chi^k_d(G)\leq c+1$.

For the third inequality, let $Q$ be a $k,d$-clique in $G$ with $|Q|=\omega^k_d(G)$ and let $f:V(G)\to[c]$ be a $k,d$-proper coloring with $c=\chi^k_d(G)$. By definition, every $k$-subset of $Q$ lies on a shortest path of $G$ of length at most $d$, and hence cannot be monochromatic. Thus, each color class of $f$ contains at most $k-1$ vertices from $Q$ and hence $|Q|\leq c(k-1)$, which implies the desired result.
\end{proof}

\begin{remark}\label{remark_perfect}
Notice that taking $k=2$ and $d=1$ in (\ref{equation_clique_bound_on_chromatic}) yields the standard clique lower bound $\omega(G)\leq\chi(G)$. Recall that a graph is \emph{perfect} if for all induced subgraphs $H$ of $G$ we have $\omega(H)=\chi(H)$. Thus, (\ref{equation_clique_bound_on_chromatic}) suggests a $k,d$-analogue of perfect graphs, which we return to in Section \ref{section_kd_clique_lower_bound}.
\end{remark}

\subsection{Subgraphs and monotonicity}

Although the chromatic number does not increase when passing to subgraphs, the opposite can be true of the $k,d$-chromatic number. For example, $P_3$ is a proper subgraph of $P_3^2\cong K_3$ with $\alpha^3_2(P_3)=2$ and $\alpha^3_2(K_3)=3$ (because there are no shortest paths of length $2$ in $K_3$). Furthermore, $\chi^3_2(P_3)=2$ whereas $\chi^3_2(P_3^2)=1$. This behavior of $\chi^k_d$ with respect to subgraphs, also exhibited by the other $k,d$-invariants discussed here, occurs because shortest paths in subgraphs need not be shortest paths in containing graphs. However, the $k,d$-invariants behave as expected when passing to isometric subgraphs. Let us note that Lemma \ref{lemma_isometric}(1)(a) and Lemma \ref{lemma_isometric}(2)(a) appeared previously in \cite[Lemma 2.6]{MR4854543}.

\begin{lemma}\label{lemma_isometric}
Suppose $k$ and $d$ are positive integers, $H$ is a subgraph of $G$ and $A\subseteq V(H)$.
\begin{enumerate}
    \item If $H$ is an isometric subgraph then the following hold.
    \begin{enumerate}
        \item If $A$ is a $k,d$-independent set in $G$ then $A$ is a $k,d$-independent set in $H$.
        \item If $A$ is a $k,d$-clique in $G$ then $A$ is a $k,d$-clique in $H$.
        \item If $A$ is a $k,d$-dominating set in $G$ then $A$ is a $k,d$-dominating set in $H$.
    \end{enumerate}
    \item If $H$ is a convex subgraph then the following hold.
    \begin{enumerate}
        \item $A$ is a $k,d$-independent set in $G$ if and only if $A$ is a $k,d$-independent set in $H$.
        \item $A$ is a $k,d$-clique in $G$ if and only if $A$ is a $k,d$-clique in $H$.
        \item $A$ is a $k,d$-dominating set in $G$ if and only if $A$ is a $k,d$-dominating set in $H$.
    \end{enumerate}
\end{enumerate}
\end{lemma}

From the previous lemma we obtain two corollaries.

\begin{corollary}\label{corollary_isometric_colorings}
Suppose $c$, $k$ and $d$ are positive integers with $d\geq k-1$. If $H$ is an isometric subgraph of $G$ and $f:V(G)\to [c]$ is a $k,d$-proper $c$-coloring of $G$ then $f\restrict V(H)$ is a $k,d$-independent $c$-coloring of $H$.
\end{corollary}
\begin{proof}
Let $H$ be an isometric subgraph of $G$ and let $f:V(G)\to[c]$ be a $k,d$-independent $c$-coloring of $G$.  
For each color $i\in[c]$, the color class $f^{-1}(i)$ is a $k,d$-independent set in $G$, and hence its restriction $(f\restrict V(H))^{-1}(i)=f^{-1}(i)\cap V(H)$ is a $k,d$-independent set in $H$ by Lemma~\ref{lemma_isometric}(1)(a).  
Therefore $f\restrict V(H)$ is a $k,d$-independent $c$-coloring of $H$, and consequently
\(\chi^k_d(H)\le\chi^k_d(G)\).
\end{proof}

\begin{corollary}\label{proposition_isometric_subgraphs}
    If $H$ is an isometric subgraph of $G$ then
    \[\alpha^k_d(H)\leq \alpha^k_d(G), \ \ \omega^k_d(H)\leq\omega^k_d(G),\ \ \chi^k_d(H)\leq\chi^k_d(G),\ \ \text{and}\ \  \gamma^k_d(H)\geq\gamma^k_d(G).\]
\end{corollary}

Next we record the relationships between particular $k,d$-invariants for different values of $k$ and $d$.

\begin{proposition}[Monotonicity in $k$ and $d$]
Suppose $G$ is a finite simple connected graph, and let $k,d,k',d'$ be positive
integers satisfying
\[
1\le k-1\le d\le \diam(G)
\qquad\text{and}\qquad
1\le k'-1\le d'\le \diam(G).
\]
If $d'\le d$ and $k'\ge k$, then the following hold:
\[\alpha^k_d(G)\leq \alpha^{k'}_{d'}(G),\ \  \chi^k_d(G)\ge \chi^{k'}_{d'}(G),\ \ \omega^k_d(G)\geq \omega^{k'}_{d'}(G),\ \ \text{and}\ \ \gamma^k_d(G)\le \gamma^{k'}_{d'}(G).
\]
\end{proposition}

\begin{proof}
    This follows from the relationship between $k,d$-versions of these kinds of special sets of vertices for various values of $k$ and $d$. For example, as in \cite[Lemma 2.6]{MR4854543}, when $d'\le d$ and $k'\ge k$, every $k,d$-independent set in $G$ is necessarily $k',d'$-independent in $G$. The remaining details of the proof are left to the reader.
\end{proof}

\subsection{Two greedy bounds}\label{section_greed_bounds}

Let us consider two different greedy bounds on $\chi^k_d(G)$. The first follows easily from known results on $k$-uniform hypergraphs, but the second exploits the underlying metric structure of the graph $G$.

\begin{corollary}[Hypergraph greedy bound]\label{corollary_greedy_hypergraph_bound}
    For positive integers $k$ and $d$ we have
    \[\chi^k_d(G)\leq\Delta(\H_{k,d}(G))+1.\]
\end{corollary}

\begin{proposition}[Metric greedy bound]\label{proposition_greedy_metric_bound}
Suppose $k$ and $d$ are positive integers with $1\leq k-1\leq d$. If $G$ is a finite graph then
\[\chi^k_d(G)\leq\Delta(G^{d-k+2})+1\leq\Delta\frac{(\Delta-1)^{d-k+2}-1}{\Delta-1}+1.\]
\end{proposition}

\begin{proof}
The second inequality is standard, so we only prove the first inequality. We will use induction on the number of vertices to prove that if $G$ is a finite graph then $G$ has a $k,d$-proper $\Delta(G^{d-k+2})+1$-coloring. Clearly the result holds when $G$ has one vertex. Suppose the result holds for graphs with $n$ vertices. Let $G$ be a graph with $n+1$ vertices such that $\Delta(G^{d-k+2})=c$. Let $H$ be the graph obtained by removing one vertex $v$ of $G$ as well as any edges containing $v$. Then $H$ has $n$ vertices and $\Delta(H^{d-k+2})\leq c$. Hence, by the inductive hypothesis, there is a $k,d$-independent $c+1$-coloring $f:V(H)\to\{1,\ldots,c+1\}$. Since $\Delta(G^{d-k+2})=c$, the set $f(N^G_{d-k+2}(v))$ contains at most $c$ of the $c+1$ available colors. We extend $f$ to a coloring $g:V(G)\to[c+1]$ by letting $g(v)$ be some color not in $f(N^G_{d-k+2}(v))$. To show that $g$ is a $k,d$-independent $c+1$-coloring of $G$, suppose $P$ is a shortest path in $G$ of length at most $d$ containing $v$. Since $f$ is a $k,d$-independent $c+1$-coloring on $H$, and since $g(v)\neq g(u)$ for $u\in N^G_{d-k+2}(v)$, it follows that the only vertices of $P$ that have the same color as $v$ must be at a distance of at least $d-k+3$ from $v$. But, since $|V(P)|\leq d+1$, there are at most $d+1-(d-k+2)=k-1$ such vertices (including $v$). Therefore, $P$ does not contain $k$ vertices of the same color.
\end{proof}

The hypergraph greedy bound and the metric greedy bound on $\chi^k_d(G)$ arise from fundamentally different obstruction mechanisms. Thus, neither bound dominates the other in general, which we see in the following example.

\begin{example}
    Let us take $k=3$ and $d=2$ through this example. First, consider $K_n$. Clearly $\chi^3_2(K_n)=1$. The metric greedy bound is $\chi^3_2(K_n)\leq\Delta(K_n)+1=n$, whereas, since $\H_{3,2}(K_n)$ has no edges, the hypergraph greedy bound is
    \[\chi^3_2(K_n)\leq\Delta(\H_{3,2}(K_n))+1=1.\]
    Thus the hypergraph greedy bound is much better in this case. However, the reason the hypergraph greedy bound is better in this case is that there are no shortest paths of length $2$ in $K_n$. Let us consider a, perhaps, more satisfying example in which there are shortest paths of length $2$, but the hypergraph greedy bound is still tighter. Suppose $n$ is even and let $G=K_n\setminus M$ where $M$ is a perfect matching of $K_n$; that is, remove $n/2$ disjoint edges from $K_n$. Then $\Delta(G)=n-2$ and $\Delta(\H_{3,2}(G))=\frac{n}{2}-1$, so the hypergraph greedy bound is still tighter.
    
    Next let us consider an example, namely the star graph $K_{1,m}$, in which the metric greedy bound is better than the hypergraph greedy bound. Of course we have $\chi^3_2(K_{1,m})=2$ (choose a color for the center vertex and then use a second color for all of the leaves). Note that $d-k+2=1$, so the graph power greedy bound becomes \[\chi^3_2(K_{1,m})\leq\Delta(K_{1,m})+1=m+1.\] On the other hand, the hyperedges in $\H_{3,2}(K_{1,m})$ are precisely those $3$-subsets of $V(K_{1,m})$ that contain the center vertex, and so the hypergraph greedy bound is 
    \[\chi^3_2(K_{1,m})\leq\Delta(\H_{3,2}(K_{1,m}))+1=\binom{m}{2}+1.\]
    Thus, the graph power greedy bound is much better than the hypergraph greedy bound in this scenario.
\end{example}

\section{$k,d$-perfect graphs}\label{section_kd_clique_lower_bound}

The following definition generalizes the classical notion of perfect graph, defined in Remark \ref{remark_perfect}.

\begin{definition}\label{definition_kd_perfect}
Let $k,d$ be positive integers with $2\le k\le d+1$, and let $G$ be a graph.  
We say that $G$ is \emph{$k,d$-perfect} if
\[
\chi^k_d(H)=\left\lceil\frac{\omega^k_d(H)}{k-1}\right\rceil
\]
for every induced subgraph $H$ of $G$.
\end{definition}

Of course a graph is $2,1$-perfect if and only if it is perfect. Furthermore, $2,d$-cliques in a graph $G$ are precisely cliques in $G^d$, and $2,d$-proper colorings of $G$ are precisely the proper colorings of $G^d$. Despite these facts, and because perfectness involves induced subgraphs, $2,d$-perfectness of a graph $G$ is not equivalent to perfectness of $G^d$. Next, we will show that if $G^d$ is perfect then $G$ is $2,d$-perfect (Proposition \ref{proposition_G^d_perfect_implies_G_is_2d_perfect}), and we provide a counterexample to the converse (Proposition \ref{proposition_counter_example}). For $k\geq 3$, some results later in the paper establish that the notion of $k,d$-perfectness provides a richer notion with features that do not appear in the classical case. For example, although we will prove that paths are $k,d$-perfect for all positive integers $k$ and $d$ (see Corollary \ref{corollary_paths_kd_perfect}), we will give a complete classification of when cycles are $k,d$-perfect (see Theorem \ref{theorem_2d_perfect_cycles} and Theorem \ref{theorem_kd_perfect_cycles}) which will show, in particular, that every cycle is $3,2$-perfect. This is in stark contrast to the classical fact that for $n\geq 4$, $C_n$ is perfect if and only if $n$ is even.

\begin{proposition}\label{proposition_G^d_perfect_implies_G_is_2d_perfect}
    Suppose $d\geq 1$ is an integer. For any graph $G$, if $G^d$ is perfect then $G$ is $2,d$-perfect.
\end{proposition}

\begin{proof}
    Suppose $G^d$ is perfect. Let $H=G[S]$ be an induced subgraph of $G$. Then $S$ induces a subgraph $H'=G^d[S]$ of $G^d$. Since $G^d$ is perfect, there is a proper coloring $f:S\to [c]$ of $H'$ with $c=\omega(H')$. Let us show that $f$ is a $2,d$-proper coloring of $H$. Suppose $u,v\in S$ and $d_H(u,v)\leq d$. Since $d_G(u,v)\leq d_H(u,v)\leq d$, it follows that $u$ and $v$ are adjacent in $G^d$ and thus also in $H'$. Hence $f(u)\neq f(v)$. Since every clique in $H'$ is a $2,d$-clique in $H$, it follows that $c\leq \omega^2_d(H)$. Since $\chi^2_d(H)\leq c\leq\omega^2_d(H)$, it follows by Proposition \ref{proposition_basic} that $\chi^2_d(H)=\omega^2_d(H)$. Hence $G$ is $2,d$-perfect.
\end{proof}

\begin{proposition}\label{proposition_counter_example}
    There is a $2,2$-perfect graph $G$ with $G^2$ not perfect.
\end{proposition}

\begin{proof}
    Let $G$ be the graph obtained by gluing a $4$-cycle to a $6$-cycle along a common edge, as in Figure \ref{fig_counterexample}(A). By inspection, one can see from Figure \ref{fig_counterexample}(A-C) that $\overline{G^2}$ contains a induced $5$-cycle, and hence $G^2$ contains an induced odd hole. It follows by the strong perfect graph theorem \cite{MR2233847} that $G^2$ is not perfect.
    
    Let us verify that $G$ is $2,2$-perfect. Suppose $H$ is an induced subgraph of $G$. First suppose $H=G$. The coloring displayed in Figure \ref{fig_counterexample}(A) shows that $\chi^2_2(G)\leq 4$, and furthermore, any $2,2$-proper coloring of $G$ must assign four distinct colors to the vertices of the $4$-cycle. Hence $\chi^2_2(G)=4$. Now, a $2,2$-clique of $G$ is precisely a clique of $G^2$, and it is clear from the diagram of $G^2$ in Figure \ref{fig_counterexample}(B) that $\omega(G^2)=4$. Hence $\chi^2_2(G)=\omega^2_2(G)$.

    Now, suppose $H$ is a proper induced subgraph of $G$ that contains a cycle. Notice that $G$ has an $8$-cycle as a subgraph, but this $8$-cycle is not an induced subgraph. So, our assumption that $H$ is a proper induced subgraph of $G$ implies that $H$ does not contain an $8$-cycle. So $H$ either contains a $4$-cycle or a $6$-cycle. Since $H$ is a subgraph of the fixed graph $G$ shown in Figure \ref{fig_counterexample}(A), it follows that if $H$ contains a $4$-cycle, then $\chi^2_2(H)=4$ and $\omega^2_2(H)=\omega(H^2)=4$. Furthermore, if $H$ contains a $6$-cycle but not a $4$-cycle, then $\chi^2_2(H)=3$ and $\omega^2_2(H)=\omega(H^2)=3$.

    Suppose $H$ does not contain a cycle. Then $H$ is a forest. Suppose $\omega(H^2)=1$. then $H$ is edgeless and so $\chi^2_2(H)=\omega^2_2(H)=1$. Suppose $\omega(H^2)=2$. Then $H$ has edges but contains no path of length $2$ as a subgraph. Hence $\chi^2_2(H)=2$. Suppose $\omega(H^2)=3$. Then $H$ contains at least one path of length $2$, and since $H$ is an induced subgraph of the fixed graph $G$, it follows by inspection that $\chi^2_2(H)=3$.\end{proof}

\begin{figure}
\begin{subfigure}{0.31\textwidth}
\centering
\begin{tikzpicture}[
every node/.style={
circle,
draw,
fill=white,
minimum size=4mm,
inner sep=0pt,
font=\scriptsize
}
]

\node[fill=CBblue] (v0) at (90:1) {0};
\node[fill=CByellow] (v1) at (90-360/6:1) {1};
\node[fill=CBgreen] (v2) at (90-2*360/6:1) {2};
\node[fill=CBblue] (v3) at (90-3*360/6:1) {3};
\node[fill=CByellow] (v4) at (90-4*360/6:1) {4};
\node[fill=CBorange] (v5) at (90-5*360/6:1) {5};

\node[fill=CBorange] (v6) at ($(v1)!1!-90:(v0)$) {6};
\node[fill=CBgreen] (v7) at ($(v0)!1!90:(v1)$) {7};

\draw (v0)--(v1);
\draw (v1)--(v2);
\draw (v2)--(v3);
\draw (v3)--(v4);
\draw (v4)--(v5);
\draw (v5)--(v0);

\draw (v0)--(v7);
\draw (v7)--(v6);
\draw (v6)--(v1);

\end{tikzpicture}\caption{$G$}
\end{subfigure}\begin{subfigure}{0.31\textwidth}
\centering\begin{tikzpicture}[
every node/.style={
circle,
draw,
fill=white,
minimum size=4mm,
inner sep=0pt,
font=\scriptsize
}
]

\foreach \i in {0,...,5}
{
    \node (v\i) at (90-360*\i/6:1) {\i};
}

\node (v6) at ($(v1)!1!-90:(v0)$) {6};
\node (v7) at ($(v0)!1!90:(v1)$) {7};

\draw (v0)--(v1);
\draw (v1)--(v2);
\draw (v2)--(v3);
\draw (v3)--(v4);
\draw (v4)--(v5);
\draw (v5)--(v0);

\draw (v0)--(v7);
\draw (v7)--(v6);
\draw (v6)--(v1);

\draw (v0)--(v2);
\draw (v1)--(v3);
\draw (v2)--(v4);
\draw (v3)--(v5);
\draw (v4)--(v0);
\draw (v5)--(v1);

\draw (v0)--(v6);
\draw (v1)--(v7);

\draw (v5) to[bend left=20] (v7);
\draw (v2) to[bend right=20] (v6);

\end{tikzpicture}\caption{$G^2$}
\end{subfigure}\begin{subfigure}{0.31\textwidth}
\centering\begin{tikzpicture}[
every node/.style={
circle,
draw,
fill=white,
minimum size=4mm,
inner sep=0pt,
font=\scriptsize
}
]

\node (v3) at (90:1.2) {3};
\node (v7) at (90-360/5:1.2) {7};
\node (v2) at (90-2*360/5:1.2) {2};
\node (v5) at (90-3*360/5:1.2) {5};
\node (v6) at (90-4*360/5:1.2) {6};

\node (v0) at (90:2) {0};

\node (v4) at (90:0.4) {4};

\node (v1) at (90:-0.4) {1};

\draw (v3)--(v7);
\draw (v7)--(v2);
\draw (v2)--(v5);
\draw (v5)--(v6);
\draw (v6)--(v3);

\draw (v0)--(v3);

\draw (v6)--(v4);
\draw (v4)--(v7);
\draw (v4)--(v1);

\end{tikzpicture}\caption{$\overline{G^2}$}
\end{subfigure}
\caption{The graph $G$, its square $G^2$, and the complement $\overline{G^2}$.}
\label{fig_counterexample}
\end{figure}

Recall that a graph $G$ is \emph{chordal} if every $n$-cycle in $G$ for $n\geq 4$ has a chord, and every chordal graph is perfect \cite{MR232694}. Although powers of chordal graphs are not chordal in general, odd powers of chordal graphs are chordal \cite{MR704427}. Thus we obtain the following from Proposition \ref{proposition_G^d_perfect_implies_G_is_2d_perfect}.
\begin{corollary}\label{corollar_chordal}
    Every chordal graph is $2,d$-perfect for all odd positive integers $d$.
\end{corollary}

Powers of strongly chordal graphs are strongly chordal \cite{MR941786, MR1206558}, where a graph $G$ is \emph{strongly chordal} if it is chordal and every $n$-cycle in $G$ of even length $\geq 6$ has an odd chord.
\begin{corollary}\label{corollary_strongly_chordal}
    Let $d\geq 1$ be an integer. Every strongly chordal graph is $2,d$-perfect; that is, $\chi^2_d(H)=\omega^2_d(H)$ for all induced subgraphs $H$ of $G$. Hence, trees, block graphs and interval graphs are all $2,d$-perfect.
\end{corollary}

\begin{proof}
    Suppose $G$ is a strongly chordal graph. Then $G^d$ is strongly chordal and thus perfect. Hence, by Proposition \ref{proposition_G^d_perfect_implies_G_is_2d_perfect}, $G$ is $2,d$-perfect.
\end{proof}

In general, perfectness and $2,d$-perfectness are not comparable properties. Perfectness does not imply $2,2$-perfectness. For example, $C_{10}$ is perfect but not $2,2$-perfect because $\omega^2_2(C_{10})=3<4=\chi^2_2(C_{10})$. Conversely $2,2$-perfectness does not imply perfectness. For example $C_5$, is $2,2$-perfect since $\omega^2_2(C_5)=5=\chi^2_2(C_5)$ and since any induced subgraph of $C_5$ is either all of $C_5$ or is a disjoint union of paths. But, of course, $C_5$ is not perfect. Similarly, for $k>3$, perfectness and $k,d$-perfectness are not comparable. For example, $C_5$ is $3,2$-perfect but not perfect.

\begin{proposition}\label{proposition_trees_are_32_perfect}
Every finite tree $T$ is $3,2$-perfect.
\end{proposition}

\begin{proof}
Let $T$ be a tree. First, let us argue that $\ceil{\frac{\omega^3_2(T)}{2}}=\chi^3_2(T)$. If $T$ has one or two vertices the result holds. Suppose $|V(T)|\geq 3$. We define a coloring $f:V(T)\to[2]$ as follows. Choose a root $r\in V(T)$ and define
\[f(v_i)=\begin{cases}
    1&\text{if $d_T(r,v_i)$ is even}\\
    2&\text{otherwise}
\end{cases}.\]
Let $P$ be any shortest path of length $2$ in $T$. Then $V(P)$ either intersects two or three distinct levels. In either case, by definition of $f$, $P$ contains two vertices of one color and one vertex of the other. Hence $f$ is a $3,2$-proper coloring and $\chi^3_2(T)\leq 2$. Since $|V(T)|\geq 3$, we have $\omega^3_2(T)=3$ and hence by Proposition \ref{proposition_basic}, we have
\[\ceil{\frac{\omega^3_2(T)}{2}}=\ceil{\frac{3}{2}}=2\leq\chi^3_2(T)=2.\]
Since every induced subgraph $H$ of $T$ is a forest, the same argument applies to each component of $H$, and thus, 
\[\ceil{\frac{\omega^3_2(H)}{2}}=\chi^3_2(H).\]
Therefore $T$ is $3,2$-perfect.
\end{proof}

\section{Coloring paths and $k,d$-perfection of paths}\label{section_coloring paths}

In order to find an exact formula for $\chi^k_d(P_n)$ we must first review some basic facts about $k,d$-independent subsets of $P_n$ from \cite{MR4854543}. Throughout the rest of the paper we assume that the vertices of $P_n$ are positive integers $V(P_n)=\{1,\ldots,n\}$ and the vertex set of the two-way infinite path graph $P_\infty$ is $V(P_\infty)=\Z$.

Let us identify the canonical largest $k,d$-independent sets in $P_n$ and the resulting exact formula for $\alpha^k_d(P_n)$. The following lemma appears in \cite{MR3985624}, but recall that the notion of $k,d$-independent set is referred to using the term $k$-general $d$-position set.

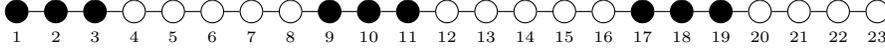
\begin{figure}[h]
\centering
\begin{tikzpicture}[
    x=0.52cm,
    every node/.style={
        circle,
        draw,
        minimum size=3mm,
        inner sep=0pt
    }
]

\node[fill=black]    (v1)  at (1,0) {};
\node[fill=black]    (v2)  at (2,0) {};
\node[fill=black]    (v3)  at (3,0) {};
\node[fill=white] (v4)  at (4,0) {};
\node[fill=white] (v5)  at (5,0) {};
\node[fill=white] (v6)  at (6,0) {};
\node[fill=white]   (v7)  at (7,0) {};
\node[fill=white]   (v8)  at (8,0) {};

\node[fill=black]    (v9)  at (9,0) {};
\node[fill=black]    (v10) at (10,0) {};
\node[fill=black]    (v11) at (11,0) {};
\node[fill=white] (v12) at (12,0) {};
\node[fill=white] (v13) at (13,0) {};
\node[fill=white] (v14) at (14,0) {};
\node[fill=white]   (v15) at (15,0) {};
\node[fill=white]   (v16) at (16,0) {};

\node[fill=black]    (v17) at (17,0) {};
\node[fill=black]    (v18) at (18,0) {};
\node[fill=black]    (v19) at (19,0) {};
\node[fill=white] (v20) at (20,0) {};
\node[fill=white] (v21) at (21,0) {};
\node[fill=white] (v22) at (22,0) {};
\node[fill=white]   (v23) at (23,0) {};

\draw (v1)--(v2)--(v3)--(v4)--(v5)--(v6)--(v7)--(v8)
      --(v9)--(v10)--(v11)--(v12)--(v13)--(v14)--(v15)--(v16)
      --(v17)--(v18)--(v19)--(v20)--(v21)--(v22)--(v23);

\node[below=5pt, font=\tiny, draw=none] at (v1)  {1};
\node[below=5pt, font=\tiny, draw=none] at (v2)  {2};
\node[below=5pt, font=\tiny, draw=none] at (v3)  {3};
\node[below=5pt, font=\tiny, draw=none] at (v4)  {4};
\node[below=5pt, font=\tiny, draw=none] at (v5)  {5};
\node[below=5pt, font=\tiny, draw=none] at (v6)  {6};
\node[below=5pt, font=\tiny, draw=none] at (v7)  {7};
\node[below=5pt, font=\tiny, draw=none] at (v8)  {8};
\node[below=5pt, font=\tiny, draw=none] at (v9)  {9};
\node[below=5pt, font=\tiny, draw=none] at (v10) {10};
\node[below=5pt, font=\tiny, draw=none] at (v11) {11};
\node[below=5pt, font=\tiny, draw=none] at (v12) {12};
\node[below=5pt, font=\tiny, draw=none] at (v13) {13};
\node[below=5pt, font=\tiny, draw=none] at (v14) {14};
\node[below=5pt, font=\tiny, draw=none] at (v15) {15};
\node[below=5pt, font=\tiny, draw=none] at (v16) {16};
\node[below=5pt, font=\tiny, draw=none] at (v17) {17};
\node[below=5pt, font=\tiny, draw=none] at (v18) {18};
\node[below=5pt, font=\tiny, draw=none] at (v19) {19};
\node[below=5pt, font=\tiny, draw=none] at (v20) {20};
\node[below=5pt, font=\tiny, draw=none] at (v21) {21};
\node[below=5pt, font=\tiny, draw=none] at (v22) {22};
\node[below=5pt, font=\tiny, draw=none] at (v23) {23};

\end{tikzpicture}
\caption{The set $S_{4,7}\cap P_{23}$ is a largest $4,7$-independent set in $P_{23}$.}\label{figure_path_coloring}
\end{figure}

\begin{lemma}[{\cite[Proof of Theorem 3.3]{MR4854543}}]\label{corollary_largest_kgdp_set_in_P_n}
Suppose $d,k$ and $n$ are positive integers with $2\leq k\leq d+1$. The set
\[S_{k,d}=\bigcup_{i\in\Z}[(d+1)i+1,(d+1)i+k-1]\]
is $k,d$-independent in $P_\infty$. Furthermore, $S_{k,d}\cap [1,n]$ is a largest $k,d$-independent set in $P_n$.
\end{lemma}

\begin{theorem}[{\cite[Theorem 3.3]{MR4854543}}]\label{theorem_gp_paths}
Suppose $d$, $k$ and $n$ are positive integers with $d\geq 1$ and $2\leq k\leq n$. Then
\[\alpha^k_d(P_n)=\begin{cases}
  n & \text{if $1\leq d\leq k-2$}\\
  (k-1)\floor{\frac{n}{d+1}}+\min(n\bmod(d+1),k-1) & \text{if $k-1\leq d$}
\end{cases}.\]
\end{theorem}

In several results that follow we will use certain canonically defined $k,d$-proper colorings on paths. For this reason, we isolate these functions now.

\begin{lemma}\label{lemma_path_coloring}
Let $a$ and $b$ be positive integers. The function $f:\Z\to\{1,\ldots,\ceil{\frac{a}{b}}\}$ defined by\footnote{See Figure \ref{figure_path_coloring} for a depiction of this function in the case where $a=8$ and $b=3$.}
\[f(i)=\floor{\frac{(i-1)\bmod a}{b}}+1\]
satisfies $f^{-1}(j)\subseteq S_{b+1,a-1}+(j-1)b$ for all $j\in\Z$ with $1\leq j\leq \ceil{\frac{a}{b}}$. Hence $f$ is a $k,d$-proper coloring of $P_\infty$.
\end{lemma}
\begin{proof}
We have $f(i)=\floor{\frac{r}{b}}+1$ where $q$ and $r$ are integers with $(i-1)=qa+r$ and $0\leq r\leq a-1$. Suppose $f(i)=j$, then $\floor{\frac{r}{b}}=j-1$ and hence $(j-1)b\leq r\leq jb-1$. Thus we have $r=(j-1)b+t$ for some $t$ with $0\leq t\leq b-1$, and so $i=qa+r+1 = qa+1+t+(j-1)b$ where $qa+1+t\in [qa+1,qa+b]\subseteq S_{b+1,a-1}$. Hence $i\in \left(S_{b+1,a-1}+(j-1)b\right)$.
\end{proof}

\begin{figure}[h]
\centering
\begin{tikzpicture}[
    x=0.52cm,
    every node/.style={
        circle,
        draw,
        minimum size=3mm,
        inner sep=0pt
    }
]

\node[fill=CBred]    (v1)  at (1,0) {};
\node[fill=CBred]    (v2)  at (2,0) {};
\node[fill=CBred]    (v3)  at (3,0) {};
\node[fill=CByellow] (v4)  at (4,0) {};
\node[fill=CByellow] (v5)  at (5,0) {};
\node[fill=CByellow] (v6)  at (6,0) {};
\node[fill=CBblue]   (v7)  at (7,0) {};
\node[fill=CBblue]   (v8)  at (8,0) {};

\node[fill=CBred]    (v9)  at (9,0) {};
\node[fill=CBred]    (v10) at (10,0) {};
\node[fill=CBred]    (v11) at (11,0) {};
\node[fill=CByellow] (v12) at (12,0) {};
\node[fill=CByellow] (v13) at (13,0) {};
\node[fill=CByellow] (v14) at (14,0) {};
\node[fill=CBblue]   (v15) at (15,0) {};
\node[fill=CBblue]   (v16) at (16,0) {};

\node[fill=CBred]    (v17) at (17,0) {};
\node[fill=CBred]    (v18) at (18,0) {};
\node[fill=CBred]    (v19) at (19,0) {};
\node[fill=CByellow] (v20) at (20,0) {};
\node[fill=CByellow] (v21) at (21,0) {};
\node[fill=CByellow] (v22) at (22,0) {};
\node[fill=CBblue]   (v23) at (23,0) {};

\draw (v1)--(v2)--(v3)--(v4)--(v5)--(v6)--(v7)--(v8)
      --(v9)--(v10)--(v11)--(v12)--(v13)--(v14)--(v15)--(v16)
      --(v17)--(v18)--(v19)--(v20)--(v21)--(v22)--(v23);

\node[below=5pt, font=\tiny, draw=none] at (v1)  {1};
\node[below=5pt, font=\tiny, draw=none] at (v2)  {2};
\node[below=5pt, font=\tiny, draw=none] at (v3)  {3};
\node[below=5pt, font=\tiny, draw=none] at (v4)  {4};
\node[below=5pt, font=\tiny, draw=none] at (v5)  {5};
\node[below=5pt, font=\tiny, draw=none] at (v6)  {6};
\node[below=5pt, font=\tiny, draw=none] at (v7)  {7};
\node[below=5pt, font=\tiny, draw=none] at (v8)  {8};
\node[below=5pt, font=\tiny, draw=none] at (v9)  {9};
\node[below=5pt, font=\tiny, draw=none] at (v10) {10};
\node[below=5pt, font=\tiny, draw=none] at (v11) {11};
\node[below=5pt, font=\tiny, draw=none] at (v12) {12};
\node[below=5pt, font=\tiny, draw=none] at (v13) {13};
\node[below=5pt, font=\tiny, draw=none] at (v14) {14};
\node[below=5pt, font=\tiny, draw=none] at (v15) {15};
\node[below=5pt, font=\tiny, draw=none] at (v16) {16};
\node[below=5pt, font=\tiny, draw=none] at (v17) {17};
\node[below=5pt, font=\tiny, draw=none] at (v18) {18};
\node[below=5pt, font=\tiny, draw=none] at (v19) {19};
\node[below=5pt, font=\tiny, draw=none] at (v20) {20};
\node[below=5pt, font=\tiny, draw=none] at (v21) {21};
\node[below=5pt, font=\tiny, draw=none] at (v22) {22};
\node[below=5pt, font=\tiny, draw=none] at (v23) {23};

\end{tikzpicture}
\caption{A \(4,7\)-proper coloring of \(P_{24}\).}\label{figure_path_coloring}
\end{figure}
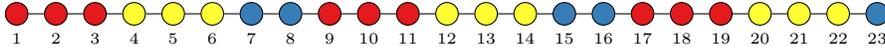

\begin{theorem}\label{theorem_chikd_paths}
Suppose $d,k$ and $n$ are positive integers with $d+1\geq k\geq 2$. Then
\[\chi^k_d(P_n)=\ceil{\frac{\min(d+1,n)}{k-1}}\]
Furthermore, \[\chi^k_d(P_\infty)=\ceil{\frac{d+1}{k-1}}.\]
\end{theorem}

\begin{proof}
We will prove the first equation in the case where $n\geq d+1$; the second equation follows by similar arguments. To see that $\chi^k_d(P_n)\leq\ceil{\frac{d+1}{k-1}}$, just notice that the function $f:V(P_n)\to\left\{1,2,\ldots,\ceil{\frac{d+1}{k-1}}\right\}$ defined by
    \[f(i)=\floor{\frac{(i-1)\bmod(d+1)}{k-1}}+1\]
is a $k,d$-proper coloring of $P_n$ by Lemma \ref{lemma_path_coloring} and Lemma \ref{lemma_isometric}(1)(a).

Let us now prove that $\chi^k_d(P_n)\geq\ceil{\frac{d+1}{k-1}}$.
Let $f:V(P_n)\to\left\{1,2,\ldots,\ceil{\frac{d+1}{k-1}}-1\right\}$ be any function. We must show that $f$ is not a proper $k,d$-proper coloring. It will suffice to show that there is an $i\in\left\{1,\ldots,\ceil{\frac{d+1}{k-1}}-1\right\}$ for which $|f^{-1}(i)\cap[1,d+1]|\geq k$. Suppose not. Then 
\begin{align}|[1,d+1]|\leq(k-1)\left(\ceil{\frac{d+1}{k-1}}-1\right).\label{equation_chi_k_d_P_n}
\end{align}
Let $q,r\in\Z$ be such that $d+1=q(k-1)+r$ where $0\leq r<k-1$. If $r=0$ then $\ceil{\frac{d+1}{k-1}}=\frac{d+1}{k-1}$, and thus by (\ref{equation_chi_k_d_P_n}) we have $|[1,d+1]|\leq d-k+2\leq d$, which is not possible. If $r>0$ then $\ceil{\frac{d+1}{k-1}}=q+1$, and hence by (\ref{equation_chi_k_d_P_n}) we have $|[1,d+1]|\leq(k-1)q=d+1-r$, and thus $|[1,d+1]|\leq d$, a contradiction. \end{proof}

\begin{corollary}\label{corollary_paths_kd_perfect}
Let $k,d,n$ be positive integers with $d+1\ge k\ge 2$. Then every finite path $P_n$ is $k,d$-perfect. Moreover, the two-way infinite path $P_\infty$ is also $k,d$-perfect.
\end{corollary}

\begin{proof}
Recall that a graph $G$ is $k,d$-perfect if every induced subgraph $H$ of $G$ satisfies
\[
\chi^k_d(H)=\ceil{\frac{\omega^k_d(H)}{k-1}}.
\] 
We only prove the result for finite paths, since the proof for infinite paths is similar. Any induced subgraph of a path $P_n$ is a disjoint union of paths, so it suffices to verify equality for paths themselves.
By Theorem~\ref{theorem_chikd_paths}, we have
\[
\chi^k_d(P_n)=\ceil{\frac{\min(d+1,n)}{k-1}},
\]
On the other hand, a largest $k,d$-clique in $P_n$ consists of a set of vertices contained in a subpath of length at most $d$, and hence has size at most $\min(d+1,n)$. Thus $\chi^k_d(P_n)=\omega^k_d(P_n)$ for all $n$, and every finite path is $k,d$-perfect.
\end{proof}

\section{Coloring cycles and a characterization of $k,d$-perfect cycles}\label{section_coloring_cycles}

In order to find an explicit formula for $\chi^k_d(C_n)$, we must first review some basic facts about some canonical $k,d$-independent sets in $C_n$. The following definition was introduced in Clough and Douthett's work in mathematical music theory \cite{CloughDouthett}, and has been studied by many authors across different disciplines \cite{Barrett,MR4875435,MR2512671,MR1401228,MR2212108,MR2567431}. The idea is that the vertices of the sets $J^r_{n,m}$ are \emph{maximally even} in the sense that they are dispersed within $C_n$ as much as possible. Indeed, the sets $J^r_{n,m}$ are minimizers of various notions of discrete electric potential energy in the same way that $m$ electrons confined to the unit circle have a unique configuration up to rotation of minimal electric potential energy (see \cite{MR4875435} for more on this perspective).

Let us note that we take the vertices of $C_n$ as being $V(C_n)=\{0,\ldots,n-1\}$ and edges $\{i,(i+1)\bmod n\}$ for each $i\in\{0,\ldots,n-1\}$. Fix vertices $u$ and $v$ in $C_n$. The \emph{clockwise distance} from $u$ to $v$ is the least nonnegative integer congruent to $(v-u)$ mod $n$, that is
\[d(u,v)=(v-u)\bmod n.\]
Recall that the \emph{distance} (or \emph{geodesic distance}) between $u$ and $v$ is the length of a shortest path from $u$ to $v$. For any set $A=\{a_0,\ldots,a_{m-1}\}$ of vertices in $C_n$ with $|A|=m\leq n$, arranged in increasing order, we define the \emph{$A$-span of $(a_i,a_k)$} to be
\[\spn_A(a_i,a_j)=(j-i)\bmod m.\]
We define two multisets: the \emph{$k$-multispectrum of clockwise distances} of $A$ is the multiset
\[\sigma^*_k(A)=\left[\,d^*(u,v)\st \spn_A(u,v)=k\,\right]\]
and the \emph{$k$-multispectrum of geodesic distances} of $A$ is the multiset
\[\sigma_k(A)=\left[\,d(u,v)\st \spn_A(u,v)=k\,\right].\]
For a multiset $X$ we define $\supp(X)$ to be the set whose elements are those of $X$.

\begin{definition}[{Clough and Douthett \cite{CloughDouthett}}]\label{definition_me} A set $A$ of vertices in $C_n$ with $|A|=m$ is \emph{maximally even} if for each integer $k$ with $1\leq k \leq m-1$ the set $\supp(\sigma^*(A))$ consists of either a single integer or two consecutive integers.
\end{definition}

\begin{definition}[{Clough and Douthett \cite{CloughDouthett}}]\label{definition_J}
Suppose $m$, $n$ and $r$ are integers with $1\leq m\leq n$ and $0\leq r\leq n-1$. We define
\[J^r_{n,m}=\left\{\floor{\frac{ni+r}{m}}\st i\in \Z \text{ and } 0\leq i < m\right\}.\]
Sets of the form $J^r_{n,m}$ are called \emph{$J$-representations}.
\end{definition}

\begin{lemma}[{Clough and Douthett \cite[Corollary 1.2]{CloughDouthett}}]\label{lemma_spec_of_J_rep_clockwise}
Suppose $k$, $m$, $n$ and $r$ are integers with $1\leq m\leq n$ and $0\leq r\leq n-1$. If $1\leq k\leq m-1$ then
\[\supp(\sigma^*_k(J^r_{n,m}))=\left\{\floor{\frac{nk}{m}},\ceil{\frac{nk}{m}}\right\}.\]
Hence $J^r_{n,m}$ is a maximally even subset of $C_n$.
\end{lemma}

Indeed, Clough and Douthett \cite{CloughDouthett} showed that a set of vertices $A$ of $C_n$ is maximally even if and only if it is of the form $J^r_{n,m}$ for some appropriate values of $m$, $n$ and $r$.

Prior to Clough and Douthett's work on maximal evenness, the following easy lemma was established by Clough and Myerson \cite[Lemma 2]{CloughAndMyerson}.

\begin{lemma}\label{lemma_sum_is_kn}
Suppose $n$ is a positive integer, $A$ is a set of vertices in $C_n$ with $|A|=m$ and $1\leq k\leq m-1$. Then
\[\sum \sigma^*_k(A)=kn.\]
\end{lemma}

The next two results were established in \cite{MR4854543}.

\begin{lemma}[{\cite[Proof of Theorem 4.6]{MR4854543}}]\label{lemma_maxeven_are_kdgenpos}
    Suppose $d,k$ and $n$ are positive integers such that $1\leq d\leq\floor{\frac{n}{2}}$ and $2\leq k\leq\floor{\frac{n}{2}}+1$. If $d\geq k-1$ and $m=\floor{\frac{(k-1)n}{d+1}}$ then $J^0_{n,m}$ is a largest $k,d$-independent set in $C_n$.
\end{lemma}

\begin{theorem}[{\cite[Theorem 4.6]{MR4854543}}]\label{theorem_alpha_kd_cycles}
    Suppose $d,k$ and $n$ are positive integers such that $1\leq d\leq\floor{\frac{n}{2}}$ and $2\leq k\leq\floor{\frac{n}{2}}+1$. Then
    \[\alpha^k_d(C_n)=\begin{cases}
        n & \text{if $1\leq d\leq k-2$}\\
        \displaystyle\floor{\frac{(k-1)n}{d+1}} & \text{if $k-1\leq d\leq\floor{\frac{n}{2}}$}
    \end{cases}.\]
\end{theorem}

Let us prove an explicit formula for $\chi^k_d(C_n)$.

\begin{theorem}\label{theorem_chikd_cycles}
    Suppose $d,k$ and $n$ are positive integers with $k\geq 2$ and $n\geq 3$. Let
    \[
    \delta=\min\left(d,\left\lfloor\frac{n}{2}\right\rfloor\right)
    \]
    and
    \[
    m=\left\lfloor\frac{(k-1)n}{\delta+1}\right\rfloor.
    \]
    Then
    \[
    \chi^k_d(C_n)=
    \begin{cases}
        1 & \text{if } 1\leq d\leq k-2,\\[0.5em]
        \displaystyle\left\lceil\frac{n}{m}\right\rceil
            & \text{if } d\geq k-1.
    \end{cases}
    \]
\end{theorem}

\begin{proof}
Without loss of generality we can assume $\delta=d\leq\floor{\frac{n}{2}}$. (In general, when $d\geq\diam(G)$, a set of vertices in $G$ is $k,d$-independent if and only if it is $k,\diam(G)$-independent. Hence, the same argument works in the case where $d>\floor{\frac{n}{2}}$, but we need to replace $d$ with $D=\diam(C_n)=\floor{\frac{n}{2}}$.)

Suppose $1\leq d\leq k-2$. In this case we have $\alpha^k_d(C_n)=n$ and hence there is a $k,d$-independent $1$-coloring of $C_n$. Thus, $\chi^k_d(C_n)=1$.

Suppose $k-1\leq d\leq \floor{\frac{n}{2}}$. Let us show that $\chi^k_d(C_n)\geq \ceil{\frac{n}{m}}.$ A largest $k,d$-independent set of vertices in $C_n$ has size $\alpha^k_d(C_n)=m$.
This implies that the largest that a color class of any $k,d$-proper coloring can be is $m$. Thus, any $k,d$-proper coloring of $C_n$ must have at least $\ceil{\frac{n}{m}}$ colors.

\begin{figure}[h]
\centering
\begin{tikzpicture}[
    every node/.style={
        circle,
        draw,
        minimum size=3mm,
        inner sep=0pt
    }
]

\def\n{16}
\def\r{1.6}

\node[fill=CBred]       (v0)  at (90:\r) {};
\node[fill=CByellow]    (v1)  at (90-360*1/\n:\r) {};
\node[fill=CBred]       (v2)  at (90-360*2/\n:\r) {};
\node[fill=CByellow]    (v3)  at (90-360*3/\n:\r) {};
\node[fill=CBblue]      (v4)  at (90-360*4/\n:\r) {};
\node[fill=CBred]       (v5)  at (90-360*5/\n:\r) {};
\node[fill=CByellow]    (v6)  at (90-360*6/\n:\r) {};
\node[fill=CBblue]       (v7)  at (90-360*7/\n:\r) {};
\node[fill=CBred]    (v8)  at (90-360*8/\n:\r) {};
\node[fill=CByellow]      (v9)  at (90-360*9/\n:\r) {};
\node[fill=CBred]       (v10) at (90-360*10/\n:\r) {};
\node[fill=CByellow]    (v11) at (90-360*11/\n:\r) {};
\node[fill=CBblue]      (v12) at (90-360*12/\n:\r) {};
\node[fill=CBred]      (v13) at (90-360*13/\n:\r) {};
\node[fill=CByellow]      (v14) at (90-360*14/\n:\r) {};
\node[fill=CBblue]      (v15) at (90-360*15/\n:\r) {};

\draw (v0)--(v1)--(v2)--(v3)--(v4)--(v5)--(v6)
      --(v7)--(v8)--(v9)--(v10)--(v11)--(v12)--(v13)--(v14)--(v15)--(v0);

\node[font=\tiny, draw=none] at (90:1.15)              {0};
\node[font=\tiny, draw=none] at (90-360*1/\n:1.15)    {1};
\node[font=\tiny, draw=none] at (90-360*2/\n:1.15)    {2};
\node[font=\tiny, draw=none] at (90-360*3/\n:1.15)    {3};
\node[font=\tiny, draw=none] at (90-360*4/\n:1.15)    {4};
\node[font=\tiny, draw=none] at (90-360*5/\n:1.15)    {5};
\node[font=\tiny, draw=none] at (90-360*6/\n:1.15)    {6};
\node[font=\tiny, draw=none] at (90-360*7/\n:1.15)    {7};
\node[font=\tiny, draw=none] at (90-360*8/\n:1.15)    {8};
\node[font=\tiny, draw=none] at (90-360*9/\n:1.15)    {9};
\node[font=\tiny, draw=none] at (90-360*10/\n:1.15)   {10};
\node[font=\tiny, draw=none] at (90-360*11/\n:1.15)   {11};
\node[font=\tiny, draw=none] at (90-360*12/\n:1.15)   {12};
\node[font=\tiny, draw=none] at (90-360*13/\n:1.15)   {13};
\node[font=\tiny, draw=none] at (90-360*14/\n:1.15)   {14};
\node[font=\tiny, draw=none] at (90-360*15/\n:1.15)   {15};

\end{tikzpicture}
\caption{A \(4,6\)-proper coloring of \(C_{16}\) using maximally even sets. The red set is $J^0_{16,6}$, the yellow set is $J^0_{16,6}+1$, and the blue vertices are contained in $J^0_{16,6}+2$.}
\label{figure_cycle_3color_asymmetric}
\end{figure}
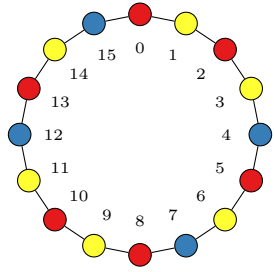

We define a $k,d$-independent $c$-coloring $f:V(C_n)\to\{1,\ldots,c\}$ where $c=\ceil{\frac{n}{m}}$ as follows. The maximally even set 
\[J^0_{n,m}=\left\{\floor{\frac{ni}{m}}\st 0\leq i<m\right\}\] 
is a $k,d$-independent set in $C_n$. The sets
\[A_R=\left\{(j+R)\bmod n\st j\in J^0_{n,m}\right\},\]
for $0\leq R<\floor{\frac{n}{m}}$, being rotations of $J^0_{n,m}$, are $k,d$-independent in $C_n$, and since $\supp(\sigma_1^*(J^0_{n,m}))=\left\{\floor{\frac{n}{m}},\ceil{\frac{n}{m}}\right\}$, they are pairwise disjoint. Let $A=\bigcup_{0\leq R<\floor{\frac{n}{m}}}A_R$. For $v\in A$ we define $f(v)=R+1$ if $v\in A_R$. If $m\mid n$, then $A=V(C_n)$ and thus $f$ is a $k,d$-independent $c$-coloring, where $c=\frac{n}{m}$. If $m\nmid n$, then $c=\floor{\frac{n}{m}}+1$, and it remains to define $f$ on $V(C_n)\setminus A$. Since $V(C_n)\setminus A$ is a subset of a rotation of the $k,d$-independent set $J^0_{n,m}$, it follows that $V(C_n)\setminus A$ is a $k,d$-independent set in $C_n$. We let $f(v)=\floor{\frac{n}{m}}+1$ for $v\in V(C_n)\setminus A$ and note that $f$ is a $k,d$-independent $c$-coloring since all of its color classes are $k,d$-independent sets by Lemma \ref{lemma_maxeven_are_kdgenpos}.
\end{proof}

We now characterize the $2,d$-perfect cycles. When $d=1$, this reduces to the classical characterization of perfect cycles, and the theorem below recovers the familiar condition that an $n$-cycle is perfect if and only if $n\le 3$ or $n$ is even.

\begin{lemma}\label{lemma_2d_clique_number_of_cycles}
Suppose $d\geq 1$ and $n\geq 3$. Then
\[
\omega^2_d(C_n)=\begin{cases}
    n & \text{if $n\leq 2d+1$}\\
    d+1 & \text{if $n\geq 2d+2$.}
\end{cases}
\]
\end{lemma}

\begin{theorem}\label{theorem_2d_perfect_cycles}
    Suppose $d\geq 1$ and $n\geq 4$. Then $C_n$ is $2,d$-perfect if and only if either,
    \begin{enumerate}
        \item $n\leq 2d+1$ or 
        \item $n\geq 2d+2$ and $(d+1)\mid n$.
    \end{enumerate}
\end{theorem}

\begin{proof}
    Every induced subgraph of $C_n$ is either a disjoint union of paths or $C_n$ itself. Since paths are $2,d$-perfect, $C_n$ is $2,d$-perfect if and only if the $2,d$-clique lower bound $\omega^2_d(C_n)$ is equal to $\chi^2_d(C_n)=\ceil{n/m(n)}$, where 
    \[m(n)=\floor{\frac{n}{d+1}}.\] If $n\leq 2d+1$ then every pair of vertices lies on a shortest path of length at most $d$ and hence $\omega^2_d(C_n)=n=\chi^2_d(C_n)$. 

    Now suppose $n\ge 2d+2$. Then $\omega^2_d(C_n)=d+1$. If $(d+1)\mid n$, then $m(n)=n/(d+1)$ and
\[
\chi^2_d(C_n)=\left\lceil \frac{n}{m(n)}\right\rceil=d+1=\omega^2_d(C_n).
\]
If $(d+1)\nmid n$, then $m(n)<n/(d+1)$, so $n/m(n)>d+1$ and hence
\[
\chi^2_d(C_n)\ge d+2>\omega^2_d(C_n).
\]
Thus $C_n$ is $2,d$-perfect if and only if $(d+1)\mid n$.
\end{proof}

Next we will characterize the $k,d$-perfect cycles for $k\ge 3$. In contrast to the case
$k=2$, where perfection is governed by a simple divisibility condition, the behavior
for larger $k$ depends on the arithmetic of the ratio $(d+1)/(k-1)$. Depending on
this arithmetic, non-perfect cycles occur either periodically or for only finitely
many values of $n$. We obtain an explicit formula for $\chi^k_d(C_n)$ (Theorem \ref{theorem_chikd_cycles}) and for $\omega^k_d(C_n)$ (in Proposition \ref{proposition_omegakd_of_C_n}), which lead to a complete classification of which cycles are $k,d$-perfect.

First, let's compute the $k,d$-clique number of $C_n$. In the next proof, we will see that there is an exceptional case where $k=3$, $n=4$ and $\omega^3_d(C_4)=4$ for any $d\geq 2$, but otherwise, a largest $k,d$-clique in $C_n$ can be found as follows. If $n\leq k-2$ then $V(C_n)$ is a largest $k,d$-clique. If the size of a largest shortest path in $C_n$ is at most $k-1$, that is $\floor{\frac{n}{2}}+1\leq k-1$, then $A_k=\{0,\ldots,k-2\}$ is a largest $k,d$-clique in $C_n$. Finally, when the set $A_k$ is a proper subset of a shortest path of $C_n$, that is when $\floor{\frac{n}{2}}+1>k-1$, the set $B_{k,d}=\left\{0,\ldots,\floor{\frac{n}{2}}+1\right\}\cap \{0,\ldots,d\}$ is a largest $k,d$-clique. This leads naturally to the following formulas.

\begin{proposition}\label{proposition_omegakd_of_C_n}
Let $3\le k\le d+1$. If $k=3$ and $n=4$, then
\[
\omega_d^3(C_4)=4.
\]
Otherwise,
\[\omega^k_d(C_n)=\begin{cases}
    n & \text{if $n\leq k-1$}\\
    k-1 & \text{if $\floor{\frac{n}{2}}+1\leq k-1$}\\
    \min\left(\floor{\frac{n}{2}}+1,d+1\right) & \text{if $\floor{\frac{n}{2}}+1>k-1$.}
\end{cases}\]
Equivalently,
\[
\omega_d^k(C_n)=
\begin{cases}
n & \text{if } n\le k-1,\\[1ex]
k-1 & \text{if } k\le n\le 2k-3,\\[1ex]
\left\lfloor\frac{n}{2}\right\rfloor+1
    & \text{if } 2k-2\le n\le 2d-1,\\[2ex]
d+1
    & \text{if } n\ge 2d.
\end{cases}
\]
\end{proposition}

\begin{proof}
If $k=3$ and $n=4$ then $V(C_4)$ is a $3,d$-clique since every triple of vertices lies on a shortest path of length $2$.

Suppose $n\leq k-1$. Then the entire vertex set $V(C_n)=\{0,\ldots,n-1\}$ has no $k$-subsets and is thus trivially a $k,d$-clique. Hence $\omega^k_d(C_n)=n$.

Suppose $\floor{\frac{n}{2}}+1\leq k-1$. Let $A_k=\{0,\ldots,k-2\}$. Since $A_k$ has no $k$-subsets, it is trivially a $k,d$-clique. We must argue that no larger set is a $k,d$-clique. Suppose $S$ is a set of vertices with $|S|\geq k$. Since $\floor{\frac{n}{2}}+1\leq k-1$, every longest shortest paths of $C_n$ contains at most $k-1$ vertices. Hence $S$ is not a $k,d$-clique. This shows that $\omega^k_d(C_n)=k-1$.

Suppose $k-1<\floor{\frac{n}{2}}+1\leq d+1$. Let $B_n=\left\{0,\ldots,\floor{\frac{n}{2}}\right\}$. Since $\floor{\frac{n}{2}}\leq d$, the entire set $B_n$ is contained in a shortest path of $C_n$ of length $d$, and so $B_n$ is a $k,d$-clique. Let us show that there is no larger $k,d$-clique. Suppose $S$ is a set of vertices with $|S|\geq\floor{\frac{n}{2}}+2$. Since $n\geq 2k-2$ we have $|S|\geq \floor{\frac{n}{2}}+2\geq k+1$. Since $k\geq 3$ we have $|S|\geq 4$. Let us first assume $|S|=4$. Since $\floor{\frac{n}{2}}+2\leq 4$ we have $n\in\{4,5\}$. Note that we can assume $n=5$ since the case $n=4$ has been handled ealier. Since $|S|=4$ and every shortest path of $C_5$ has length $2$, it is easy to see that $S$ contains a triple that does not live on any shortest path of $C_5$. Let us assume $|S|\geq 5$. For each vertex $v\in S$ there is a unique $v_{cw}\in C_n$ that is the endpoint of the clockwise shortest path of length $\floor{\frac{n}{2}}$ starting at $v$. So, there must be a pair $\{x,y\}\subseteq S$ such that $d_{C_n}(x,y)=\floor{\frac{n}{2}}$, because otherwise, for each $v\in S$ the vertex $v_{cw}$ is not in $S$ and the map $v\mapsto v_{cw}$ is injective, hence $C_n$ would have at least $2|S|\geq\left(\floor{\frac{n}{2}}+2\right)$ vertices, which is not possible. Let $P$ be a shortest path in $C_n$ from $x$ to $y$. Since $|S-\{x,y\}|\geq\floor{\frac{n}{2}}$ and $|V(P)|-\{x,y\}|=\floor{\frac{n}{2}}-1$, it follows by a similar argument that there is a second pair $\{s,t\}\subseteq S-\{x,y\}$ with $d_{C_n}(s,t)=\floor{\frac{n}{2}}$ (otherwise $C_n$ would have more than $n$ vertices). Now, since $|S|\geq 5$ we may choose $a\in S-\{x,y,s,t\}$. Then $a$ lives in a cyclic interval with endpoints among $\{x,y,s,t\}$. Without loss of generality, say $a$ is in the cylic interval $(y,s)_{C_n}$, then we choose $b=x$ and $c=t$ to be the other members of the corresponding diametral pairs. Since $b$ and $c$ form a diametral pair, every $b,c$-path has length at least $\floor{\frac{n}{2}}$, and because $a$ lies outside every shortest $b,c$-path, any path containing $\{a,b,c\}$ must properly contain a shortest $b,c$-path. Hence, any path in $C_n$ containing $\{a,b,c\}$ is too long to be a shortest path. Since $\{a,b,c\}\subseteq S$ and $k\geq 3$, this shows that $S$ is not a $k,d$-clique in $C_n$. Hence $\omega^k_d(C_n)=\floor{\frac{n}{2}}+1$.

Suppose $k-1<\floor{\frac{n}{2}}+1$ and $d+1<\floor{\frac{n}{2}}+1$. Then $d<\diam(C_n)$, and so, every pair of vertices at distance $d$ is joined by a unique shortest path. Let $S$ be a $k,d$-clique of $C_n$. If $|S|< k$, then $|S|\leq k-1\leq d+1,$
so assume $|S|\geq k$. Since every $k$-subset of $S$ is contained in a shortest path of length at most $d$, every pair of vertices of $S$ has distance at most $d$. Hence the diameter of $S$ is at most $d$. Choose vertices $s,t\in S$ with $d_{C_n}(s,t)=\operatorname{diam}(S).$ We have $d_{C_n}(s,t)\leq d<\floor{\frac{n}{2}}$, so there is a unique shortest $s,t$ path $P$. We claim that $S\subseteq V(P)$. Suppose, to the contrary, that some vertex $z\in S$ does not lie on $P$. Then
it follows that either $d_{C_n}(s,z)>d_{C_n}(s,t)$ or $d_{C_n}(t,z)>d_{C_n}(s,t)$, contradicting the maximality of $d_{C_n}(s,t)$. Thus $S$ is contained in the shortest path $P$, which has $d+1$ vertices. Hence $|S|\leq d+1$. Since the vertex set of any shortest path of length $d$ is itself a $k,d$-clique, we conclude that $\omega_d^k(C_n)=d+1.$ Therfore, $\omega^k_d(C_n)=d+1$.
\end{proof}

Let us observe that computation of $\omega^k_d(C_n)$ simplifies nicely if one assumes that $k$ vertices can fit into a shortest path of length $d$ and that shortest paths of length $d$ exist in $C_n$.

\begin{corollary}
    If $2\leq k\leq d+1$ and $d\leq\diam(C_n)$ then $\omega^k_d(C_n)=d+1$.
\end{corollary}

In order to classify $k,d$-perfect cycles for a particular parameter regime, we will need to simplify the size of a largest $k,d$-independent set (see Theorem \ref{theorem_alpha_kd_cycles}), which the following lemma will allow us to do.

\begin{lemma}
Suppose $2\leq k$ and $n\geq 2k-2$. Then
\[
\left\lfloor
\frac{(k-1)n}{\lfloor n/2\rfloor+1}
\right\rfloor
=
\begin{cases}
2k-4,&\text{if $n$ is even and } n < 4k-6,\\
2k-3,&\text{$n$ is odd or $n\geq 4k-6$.}
\end{cases}
\]
\end{lemma}

\begin{proof}
First suppose $n=2r$, where $r\in \Z$. Then
\[
\left\lfloor
\frac{(k-1)n}{\lfloor n/2\rfloor+1}
\right\rfloor
=
\left\lfloor
\frac{2(k-1)r}{r+1}
\right\rfloor
=
\left\lfloor
2k-2-\frac{2k-2}{r+1}
\right\rfloor.
\]
Since $n\geq 2k-2$, it follows that $0<\frac{2k-2}{r+1}<2$. This implies that the floor is either $2k-3$ or $2k-4$. Furthermore, the floor is equal to $2k-4$ if and only if $1<\frac{2k-2}{r+1}<2$, which is equivalent to $n<4k-6$. This establishes that the results holds in the case where $n$ is even.

Now suppose $n=2r+1$, where $r\in\Z$. Then
\[
\left\lfloor
\frac{(k-1)n}{\lfloor n/2\rfloor+1}
\right\rfloor
=
\left\lfloor
\frac{(k-1)(2r+1)}{r+1}
\right\rfloor
=
\left\lfloor
2k-2-\frac{k-1}{r+1}
\right\rfloor.
\]
Since $n\ge2k-2$, we again have $r+1\ge k$, and therefore
\[
0<\frac{k-1}{r+1}<1.
\]
Hence the floor is $2k-3$ whenever $n$ is odd.
\end{proof}

Let us now give the classification of $k,d$-perfect cycles.

\begin{theorem}\label{theorem_kd_perfect_cycles}
Let $2\leq k\leq d+1$. 
\begin{enumerate}
\item Suppose $\diam(C_n)<d$.
    \begin{enumerate}
        \item\label{item_diam<d_1} If $\diam(C_n)\leq k-2$ then $C_n$ is $k,d$-perfect.
        \item\label{item_diam<d_2} If $\diam(C_n)>k-2$, then $C_n$ is $k,d$-perfect if and only if $n\le\left(2k-3\right)
\left\lceil
\frac{\lfloor n/2\rfloor+1}{k-1}
\right\rceil$.

    \end{enumerate}

\item\label{item_diam>=d} Suppose $d\leq\diam(C_n)$ and define \[q = \frac{d+1}{k-1}.\]

\begin{enumerate}
\item\label{item_q_in_Z} Suppose $q\in\mathbb{Z}$. Then
\[
C_n \text{ is $k,d$-perfect } \iff q \mid n.
\]
\item\label{item_q_not_in_Z} Suppose $q\notin\mathbb{Z}$. Then $C_n$ is $k,d$-perfect unless there is an integer $j$ with 
\[
1 \le j < \frac{q-1}{\lceil q\rceil - q}.
\]
such that 
\[j\ceil{q}<n<(j+1)q.\]
\end{enumerate}
\end{enumerate}
\end{theorem}

        


\begin{proof}
Since induced subgraphs of cycles are disjoint unions of paths, it follows that $C_n$ is $k,d$-perfect if and only if the clique lower bound $\ceil{\frac{\omega^k_d(C_n)}{k-1}}\leq\chi^k_d(C_n)$ is an equality.

Let us prove (\ref{item_diam<d_1}). Suppose $\diam(C_n)<d$. Additionally assume $\diam(C_n)\leq k-2$. Then $\omega^k_d(C_n)=k-1$ and $\chi^k_d(C_n)=\ceil{\frac{n}{m}}$ where $m=\floor{\frac{(k-1)n}{\delta+1}}$ and $\delta=\min\left(d,\floor{\frac{n}{2}}\right)=\floor{\frac{n}{2}}$. The $k,d$-clique lower bound is
\[1=\ceil{\frac{\omega^k_d(C_n)}{k-1}}\leq \chi^k_d(C_n)=\ceil{\frac{n}{m}}.\]
Since $\floor{\frac{n}{2}}+1\leq k-1$ by assumption, it follows that $m\geq n$ and hence $\chi^k_d(C_n)=\ceil{\frac{n}{m}}=1$. This implies that $C_n$ is $k,d$-perfect.


Let us prove (\ref{item_diam<d_2}). Suppose $\diam(C_n)<d$ and $\diam(C_n)>k-2$. It follows that $2k-2\leq n\leq 2d-1$. Furthermore, by Proposition \ref{proposition_omegakd_of_C_n} we have $\omega^k_d(C_n)=\min(\floor{\frac{n}{2}}+1,d+1)=\floor{\frac{n}{2}}+1<d+1$ and by Theorem \ref{theorem_chikd_cycles}, $\chi^k_d(C_n)=\ceil{\frac{n}{m}}$ where $m=\floor{\frac{(k-1)n}{\delta+1}}$ and $\delta=\min(d,\floor{\frac{n}{2}})=\floor{\frac{n}{2}}$. Let
\[
q(n)=\frac{\lfloor n/2\rfloor+1}{k-1}\]
and note that $m=\left\lfloor\frac{n}{q(n)}\right\rfloor$.
So, the $k,d$-clique lower bound is
\[\left\lceil\frac{\omega_d^k(C_n)}{k-1}\right\rceil
=\lceil q(n)\rceil\leq \ceil{\frac{n}{m}}=\chi^k_d(C_n).\]
Since $m\le n/q(n)$, we have $n/m\ge q(n)$, and hence $C_n$ is $k,d$-perfect if and only if 
\begin{align}
n\le m\lceil q(n)\rceil.
\label{equation_perfect_proof}
\end{align}
By the lemma, either $m=2k-3$, or $n$ is even and $n<4k-6$. In the former case, we see that
\[
C_n\text{ is $k,d$-perfect}
\iff
n\le(2k-3)\lceil q(n)\rceil.
\]
In the latter case, $2k-2\le n<4k-6$ and $n$ is even, so $1<q(n)<2$, whence $\lceil q(n)\rceil=2$. Moreover, $m=2k-4$, and therefore (\ref{equation_perfect_proof}) holds since $n\le4k-8$. Thus the same characterization holds in this case as well.

Let us consider (\ref{item_diam>=d}). Suppose $d\leq\diam(C_n)$. Let $m=\floor{\frac{n}{q}}$. We have $\chi^k_d(C_n)=\ceil{\frac{n}{m}}$ and $\omega^k_d(C_n)=d+1$ (since $d\leq\diam(C_n)$). So the $k,d$-clique lower bound is 
\[\ceil{q}\leq\ceil{\frac{n}{m}}.\]
Hence, for all relevant values of $n$, $C_n$ is not $k,d$-perfect if and only if $n>m\ceil{q}$. Consider the partition of the nonnegative reals into intervals $I_j=[jq,(j+1)q)$ where $j$ ranges over nonnegative integers. Notice that $m=j$ if and only if $n\in I_j$. Hence, if $n\in I_j$ then $C_n$ is not $k,d$-perfect if and only if $n>j\ceil{q}$.

Let us prove (\ref{item_q_in_Z}). Suppose $q\in \Z$. Fix a nonnegative integer $j$. If $n=jq$ then $C_n$ is perfect, and if $n\in(jq,(j+1)q)$ then $C_n$ is not perfect. Thus, we obtain (\ref{item_q_in_Z}).

For (\ref{item_q_not_in_Z}), suppose $q\notin \Z$. For each $j$ and for $n\in I_j$, we have that $C_n$ is not $k,d$-perfect if and only if $n\in(j\ceil{q},(j+1)q)$. Notice that the interval $(j\ceil{q},(j+1)\ceil{q})$ contains integers if and only if $j\ceil{q}+1<(j+1)q$, which is equivalent to $j<\frac{q-1}{\ceil{q}-q}$.
\end{proof}

By unraveling what Theorem \ref{theorem_kd_perfect_cycles} says in the case where $k=3$ we obtain the following.

\begin{corollary}\label{corollary_3d_perfect_cycles}
The $3,d$-perfection of cycles is described in Table~\ref{table_3d_cycles}.
\end{corollary}

\begin{table}[h]
\centering
\begin{tabular}{c|l}
\(d\) & The $n\geq 3$ for which $C_n$ is $3,d$-perfect \\ \hline
2 & all $n$ \\
3 & ($n<7$) or ($n\geq 7$ and $2\mid n$) \\4 & $n\notin\{7\}$ \\
5 & ($n < 10$ and $n\notin\{7\}$) or ($n\geq 10$ and $3\mid n$) \\
6 & $n\notin\{7,10,11,13,17\}$ \\
7 & ($n<14$ and $n\notin\{7,10,11,13\}$) or ($n\geq 14$ and $4\mid n$) \\
8 & $n\notin \{7, 10, 11, 13, 14, 15, 16, 17, 21, 22, 26, 31\}$ \\
9 & ($n<18$ and $n\notin \{7, 10, 11, 13, 14, 15, 16, 17\}$) or ($n\geq 18$ and $5\mid n$) \\
\end{tabular}
\caption{Classification of \(3,d\)-perfect cycles (See Theorem \ref{theorem_kd_perfect_cycles} and Corollary \ref{corollary_3d_perfect_cycles}).}
\label{table_3d_cycles}
\end{table}

\section{Open Questions and Future Directions}\label{section_questions}

The results of this paper suggest that even for very structured graph families, the behavior of $k,d$-invariants can be subtle and highly sensitive to both parameters. While closed formulas for $\alpha^k_d$, $\chi^k_d$ and $\omega^k_d$ are now available for paths (Section \ref{section_coloring paths}), cycles (Section \ref{section_coloring_cycles}), and powers of paths (see \cite{CodyDetore_powers}), comparatively little is known for higher-dimensional or product graphs. In particular, the interaction between geodesic structure and higher-order constraints becomes significantly more complex in graphs that admit many shortest paths between vertex pairs. This motivates the systematic study of $k,d$-invariants for classical graph families beyond trees and cycles.

\begin{question}
Compute or estimate $k,d$-invariants for other classical graph families, such as:
\begin{enumerate}
    \item hypercubes $Q_n$,
    \item complete $\ell$-partite graphs $K_{n_1,\ldots,n_\ell}$,
    \item Kneser graphs, and
    \item distance-regular graphs (including Johnson and Hamming graphs).
\end{enumerate}
Which of these families is $k,d$-perfect for various values of $k$ and $d$?
\end{question}

\begin{question}
Determine the $k,d$-independence number $\alpha^k_d(G)$ and the $k,d$-chromatic number $\chi^k_d(G)$ for Cartesian products of paths. In particular:
\begin{enumerate}
    \item What is $\chi^k_d(P_\infty\Box P_\infty)$?
    \item Find exact formulas or sharp bounds for $\chi^k_d(P_n \Box P_m)$ and $\chi^k_d(C_n \Box C_m)$.
    \item More generally, determine $\chi^k_d(P_{n_1} \Box \dots \Box P_{n_r})$ for fixed $k$ and $d$.
\end{enumerate}
\end{question}

We note that the case of $\chi^k_d(P_n \Box P_2)$, though likely difficult, should be approachable given that a formula for $\alpha^k_d(P_n \Box P_2)$ is known \cite{MR4854543}. The canonical largest $k,d$-independent (or $k$-general $d$-position) subsets of $P_n \Box P_2$ used in \cite{MR4854543} serve as natural candidates for color classes in optimal colorings.

In classical coloring theory, Brooks' Theorem and its subsequent refinements \cite{MR4688000} provide fundamental bounds on the chromatic number in terms of the maximum degree $\Delta(G)$. It is natural to investigate whether the $k,d$-chromatic number obeys similar bounds that improve the greedy bounds obtained in Section \ref{section_basic_properties}. An analogue of Brooks' Theorem is obtained for $\chi^2_d$ in \cite{MR4688000}.

\begin{question}
Fix $k$ and $d$. What is the maximum possible value of $\chi^k_d(G)$ among graphs $G$ of order $n$ and bounded maximum degree $\Delta(G)$? 
\end{question}

While the study of $k,d$-independent sets in cycles involves $1$-dimensional metric constraints, higher powers of cycles introduce have a denser edge structure that can significantly compress the available $k,d$-independent sets. We propose the following problem to explore this density transition.

\begin{openproblem}
For positive integers $k,d,\ell$, determine an explicit formula for the $k,d$-independence number $\alpha^k_d(C_n^\ell)$ and the $k,d$-chromatic number $\chi^k_d(C_n^\ell)$ of the $\ell$-th power of the cycle $C_n$, or give sharp asymptotic bounds as $n \to \infty$.
\end{openproblem}

The case $\ell=1$ is handled in Theorem~\ref{theorem_chikd_cycles}. For powers of cycles, however, the interaction between the metric constraints, the power parameter $\ell$, and the parity of $n$ appears substantially more subtle than for paths. This leads to a natural question regarding perfection in these families.

\begin{question}
Fix integers $k \ge 2$ and $d \ge 1$. For which values of $n$ and $\ell$ are the graphs $P_n^\ell$ and $C_n^\ell$ $k,d$-perfect?
\end{question}

In Proposition~\ref{proposition_trees_are_32_perfect} we established that finite trees are $3,2$-perfect. Furthermore, Theorem~\ref{theorem_2d_perfect_cycles} establishes that all cycles $C_n$ are $3,2$-perfect for $n \ge 1$. This stands in stark contrast to the classical case ($k=2, d=1$) where odd holes are the primary obstruction to perfection. This unexpected ``good behavior'' of cycles suggests that the class of $3,2$-perfect graphs may be significantly larger than the class of classical perfect graphs.

\begin{question}[The 3,2-Perfection Problem]
Which graphs are $3,2$-perfect?
\end{question}

\begin{question}
Are all finite trees $k,d$-perfect for all positive integers $k,d$? 
\end{question}

The perfect graph theorem \cite{MR309780} states that an undirected graph is perfect if and only if its complement graph is also perfect. The strong perfect graph theorem \cite{MR2233847} characterizes classical perfection via forbidden induced subgraphs (odd holes and antiholes). Our results for $k \ge 3$ show that obstructions to $k,d$-perfection in cycles, for some parameter regimes, are governed by arithmetic conditions on the ratio $(d+1)/(k-1)$.
\begin{question}
Is there an analogue of either the perfect graph theorem or the strong perfect graph theorem for $k,d$-perfection? Specifically, is the $k,d$-perfection of $G$ equivalent to some notion of perfection of the uniform complement of the hypergraph $\H_{k,d}(G)$ (see Remark \ref{remark_hypergraph})? Can the class of $k,d$-perfect graphs be characterized by excluding subgraphs?
\end{question}

\section{Acknowledgments}

The second author was partially supported by an Undergraduate Research Opportunities Program award from Virginia Commonwealth University.


\end{document}